\documentclass{ifacconf}

\newtheorem{remark}{Remark}[section]

\usepackage{graphicx}      % include this line if your document contains figures
\usepackage{natbib}
\usepackage[latin1]{inputenc}
\usepackage[english]{babel}
\usepackage{graphicx}
\usepackage{amsmath}
\usepackage{amssymb}
\usepackage{color}
\usepackage{subfigure}
\usepackage{epsfig}
\usepackage{epstopdf}
\begin{document}
\begin{frontmatter}

\title{Freeway Ramp Metering Control \\ Made Easy and Efficient%\thanksref{xxx}
}

\author[First]{Hassane Aboua\"{\i}ssa}
\author[Second,Fifth]{Michel Fliess}
\author[Third]{Violina Iordanova}
\author[Fourth,Fifth]{C\'{e}dric Join}

\address[First]{Univ. Lille Nord France, 59000 Lille, France
\\U-Artois, LGI2A (EA 3926),
Technoparc Futura \\ 62400 B\'ethune, France (e-mail:
hassane.abouaissa@univ-artois.fr)}
\address[Second]{LIX (CNRS, UMR 7161),
\'Ecole polytechnique \\ 91228 Palaiseau, France (e-mail:
Michel.Fliess@polytechnique.edu)}
\address[Third]{CETE Nord-Picardie, D\'epartement Transport Mobilit\'es \\
Groupe Syst\`emes de Transport,\\ 2 rue de Bruxelles, BP 275, 59019
Lille, France
\\ (e-mail: Violina.Iordanova@developpement-durable.gouv.fr)}
\address[Fourth]{INRIA -- Non-A \& CRAN (CNRS, UMR 7039) \\
Universit\'{e} de Lorraine, BP 239 \\ 54506 Vand\oe{}uvre-l\`es-Nancy,
France
\\ (e-mail: Cedric.Join@univ-lorraine.fr) }
\address[Fifth]{AL.I.E.N., 24-30 rue
Lionnois, BP 60120, 54003 Nancy, France \\
www.alien-sas.com}

\begin{abstract}:                % Abstract of not more than 250 words.
``Model-free'' control and the related ``intelligent''
proportional-integral (PI) controllers are successfully applied to
freeway ramp metering control. Implementing moreover the
corresponding control strategy is straightforward. Numerical
simulations on the other hand need the identification of quite
complex quantities like the free flow speed and the critical
density. This is achieved thanks to new estimation techniques where
the differentiation of noisy signals plays a key r\^{o}le. Several
excellent computer simulations are provided and analyzed.
\end{abstract}


\end{frontmatter}
%===============================================================================

\section{Introduction}
This plenary lecture aims at presenting in a clear and unified
manner recent advances due to the same authors (\cite{douai,agadir})
on two important subjects in intelligent transportation systems
(see, \textit{e.g.}, \cite{gl}, \cite{ko}, \cite{mammar}, and the
references therein):
\begin{enumerate}
\item the control of freeway ramp metering,
\item the estimation of the free-flow speed and of the critical
density.
\end{enumerate}
Freeway ramp metering control, which should alleviate congestions,
is achieved via \emph{model-free control} (\cite{esta,malo}). It
yields an \emph{intelligent proportional-integral}, or \emph{iPI},
controller which
\begin{itemize}
\item regulates the traffic flow in a most efficient way,
\item is robust with respect to quite strong disturbances,
\item is easy to tune and to implement,
\item does not need any precise mathematical modeling.
\end{itemize}
\begin{remark}Model-free control, although quite
new, has already been successfully employed in many concrete
situations:

\cite{acc,cifa10,buda,cifa12,nolcos,brest,psa,edf,michel,vil1,vil2,mounier}.\end{remark}
Computer experiments show that our control strategy behaves better
than \emph{ALINEA},\footnote{\emph{ALINEA} is an acronym of
\textit{\underline{A}sservissement \underline{LIN}\'{e}aire
d'\underline{E}ntr\'{e}e \underline{A}utorouti\`{e}re}.} which was a most
remarkable breakthrough when introduced more than twenty years ago
(\cite{hist1,hist2,hist3}).\footnote{See, \textit{e.g.},
(\cite{hist4,hist5}) for recent developments.} Despite a huge
academic literature, which utilizes most of the existing methods of
modern control theory, whether with lumped or with distributed
parameter systems, ALINEA, which is exploited in France and in many
other countries, remains until today to the best of our knowledge
the only feedback-control law for ramp metering that has been
implemented in practice.

Computer simulations, on the other hand, need some kind of precise
mathematical macroscopic modeling. They become therefore more subtle
and complex than model-free control. This necessity yields a
dichotomy which is analyzed in this paper for the first time. We
utilize here ordinary differential equations, \textit{i.e.}, a
macroscopic model of order two, due to \cite{Payne} and improved by
\cite{papageorgiou}. The corresponding model properties are quite
sensitive to parameter variations and uncertainties. The
\emph{free-flow speed} and the \emph{critical density} are estimated
here via the \emph{fundamental diagram} due to \cite{May} thanks to
recent differentiation techniques of noisy signals
(\cite{nl,mboup}). Most of the existing methods for achieving
real-time estimation employ in one way or the other the Kalman
filtering (see, \textit{e.g.}, \cite{ekf,Wang05,Wang08}). Their
computational burden seems however quite higher than ours.

Our paper is organized as follows. Model-free control and
intelligent PI controllers are presented in Section
\ref{mfc}.\footnote{See \cite{marseille} for a complete
presentation.} Section \ref{Frm} studies the application to a
concrete example of an isolated ramp metering. After reviewing the
identification techniques which are connected to the fundamental
diagram, important parameters corresponding to the same freeway are
estimated in Section \ref{pe}. Convincing computer simulations are
also analyzed in Sections \ref{Frm} and \ref{pe}. Some concluding
remarks are discussed in Section \ref{conclusion}.

\section{Model-free control: A short review}\label{mfc}
\subsection{Basics}
We restrict ourselves for simplicity's sake to a SISO system
$\mathfrak{S}$, with a single input $u$ and a single output $y$. We
do not know any global mathematical description of $\mathfrak{S}$.
We replace it by a ``phenomenological'' model, which is
\begin{itemize}
\item valid during a short time lapse,
\item said to be \emph{ultra-local},
\end{itemize}
\begin{equation}\label{F}
\boxed{y^{(\nu)} = F + \alpha u}
\end{equation}
where
\begin{itemize}
\item the differentiation order $\nu$ of $y$, which is
\begin{itemize}
\item chosen by the practitioner,
\item generally equal to $1$,
\end{itemize}
has no connection with the unknown differentiation order of $y$ in
$\mathfrak{S}$;

\item the constant parameter $\alpha$ has no \textit{a priori} precise numerical
value. It is determined by the practitioner in such a way that the
numerical values of $\alpha u$ and $y^{(\nu)}$ are of equivalent
magnitude;

\item $F$, which contains all the ``structural'' information,
depends on all the system variables including the perturbations.
\end{itemize}
%L'estimation en temps r\'{e}el de la valeur num\'{e}rique de $F$, trait\'{e}e au
%{\S} \ref{oe}, permet de r\'{e}actualiser \eqref{F} \`{a} chaque instant.

\subsection{Intelligent PI controllers}
Assume that we have a ``good'' estimate\footnote{See Section
\ref{estF}.} $[F]_e$ of $F$ and, for simplicity's sake, that $\nu =
1$ in Equation \eqref{F}.\footnote{$\nu = 1$ is an appropriate
choice for most of the concrete examples. See \cite{marseille} for
an explanation.} The desired behavior is obtained via an
\emph{intelligent proportional-integral}, or \emph{iPI}, controller
\begin{equation}\label{ipi}
\boxed{u =  - \frac{[F]_e - \dot{y}^\star + K_P e + K_I \int
e}{\alpha}}
\end{equation}
where
\begin{itemize}
\item $y^\star$ is the output reference trajectory,
\item $e = y - y^\star$ is the tracking error,
\item $K_P$, $K_I$ are the usual gains.
\end{itemize}
If $K_I = 0$, we have an \emph{intelligent proportional}, or
\emph{iP}, controller:
\begin{equation}\label{iP}
\boxed{u = - \frac{[F]_e - \dot{y}^\ast + K_P e }{\alpha}}
\end{equation}

\begin{remark}
Contrary to the situations with classic PI controllers, controllers
\eqref{ipi} and \eqref{iP} are easy to tune: they stabilize a pure
integrator.
\end{remark}

\begin{remark}
See \cite{pid} for the explanation of the strange ubiquity of
classic PIDs via the above viewpoint.\footnote{See also
\cite{marseille}.}
\end{remark}

\subsection{Estimation of $F$}\label{estF}
\subsubsection{Estimation of $\dot{y}$}
If $\nu = 1$ in Equation \eqref{F}, $[F]_e$ may be obtained via the
estimate of $\dot{y}$. Elementary differentiation filters do suffice
in this situation where the sampling is rather crude.

\subsubsection{Another technique}
Rewrite Equation \eqref{ipi} as
\begin{equation*}\label{est1}
F = - \alpha u + \dot{y}^\star - K_P e - K_I \int e
\end{equation*}
Corrupting noises are attenuated by integrating both sides on a
short time interval.\footnote{See \cite{bruit} for a mathematical
explanation.} It yields:
\begin{equation}
F_{\text{\tiny approx}} = \frac{1}{\delta}\int_{T-\delta}^{T}\left
(-\alpha u + \dot{y}^\star - K_P e -K_I \int e\right ) d\tau
\label{estimF}
\end{equation}
where $F_{\text{\tiny approx}}$ is a piecewise constant
approximation of $F$. Equation \eqref{estimF} may be easily
implemented as a discrete linear filter.

\section{Freeway ramp metering principle}\label{Frm}
\subsection{Generalities}
\begin{figure}
\begin{center}
\includegraphics[width=8.6cm]{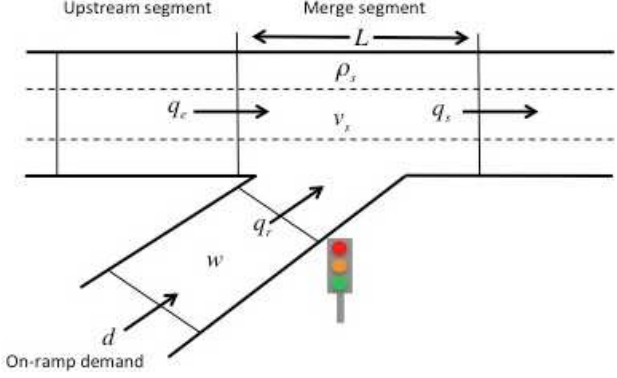}    % The printed column width is 8.4 cm.
\caption{Freeway ramp metering principle}
\label{fig:rampe}
\end{center}
\end{figure}

Consider the simple example of the freeway section depicted in Fig.
\ref{fig:rampe}:
\begin{itemize}
\item $q_r$, in $veh/h$, is the ramp flow related to the control
variable $r \in \{r_{min}, r_{max}\}$,\footnote{Just as in \cite{Hegyi_r}, we set $r=1$ for an unmetered on-ramp.} by $q_r=r\hat q_r$, where $\hat q_r=\min
\left(d+\frac{w}{T_s},Q_{sat}\min
\left(r,\frac{\rho_{max}-\rho_s}{\rho_s-\rho_c}\right)\right)$ is
the flow; \item $w$ represents the queue length in vehicles,
\item $Q_{sat}$, is the on-ramp capacity in $veh/h$,
\item $\rho_{max}$, $\rho_{c}$ are respectively the maximum and the critical density.
\end{itemize}
The \emph{ramp metering}, or \emph{admissible control}, consists  to
act on the traffic demand at the on-ramp origin in order to maintain
the traffic flow in the mainstream section close to the critical
density.\footnote{The traffic demand is assumed to be independent of
any control actions (see, \textit{e.g.},
\cite{papageorgiou,kostialos1}).}

Ramp metering strategies may be local (isolated ramp metering) or
coordinated (\cite{hist5}). Isolated ramp metering makes use of
real-time traffic measurements in the vicinity of each controlled
on-ramp in order to calculate the corresponding suitable ramp
metering flows. Coordinated ramp metering exploits the all available
measurements of the considered portions of controlled freeway. We
focus here on isolated ramp metering.

\subsection{Model-free ramp metering}
For the studied freeway section (see Fig. \ref{fig:rampe}), Equation
\eqref{F} becomes\footnote{The traffic occupancy measurements are
utilized for practical implementation.}
\begin{equation}\label{densite}
\boxed{\dot{\rho}_s(t)=F(t)-\alpha r(t)}
\end{equation}
The control variable $r(t)$ is given {\it via} the intelligent
controller iPI \eqref{ipi}:
\begin{equation}\label{cde}
\boxed{r(t)=\frac{1}{\alpha}\left[ -[F]_e +\dot \rho^\star+ K_P e +
K_I\int e \right]}
\end{equation}
where
\begin{itemize}
\item [-] $\rho^\star$ is the reference trajectory.
\item [-] $e=\rho_s-\rho^\star$ is the tracking error.
\end{itemize}
The estimation of $F$ is provided thanks to the following expression:
\begin{equation*}
\boxed{[F(k)]_e = \left[\dot{\rho}_s(k)]_e-\alpha r(k-1)\right]}
\end{equation*}
where
\begin{itemize}
\item [-]$k$ is the sampled time,
\item [-] $[\bullet]_e$ indicates an estimate of
$\bullet$.\footnote{See Section \ref{estF}.}
\end{itemize}
%Let us emphasize that, for the real-time implementation, there is no
%need of the use of the derivative technique explained in Section
%\ref{ns}. Such remark allows us to underline the following fact:
%Using the model-free  approach the control problem becomes trivial
%since the conceptual framework of the control becomes simple.
%Nevertheless, for simulation purposes, it is impossible to pass off
%of a more realistic model.

%\begin{remark}
%In contrast with the other strategies, no need here of macroscopic
%models\footnote{Recall that such models, based on three macroscopic
%variables (traffic density, mean speed and flow), are the most
%appropriate for freeway network description (See e.g.
%\cite{papageorgiou98}).}, which are based on partial differential
%equations.  The iPI \eqref{cde} allows, effectively to take into
%account the traffic complexity thanks to the estimated term $[F]_e$.
%\end{remark}
\subsection{Implementation issues}
\begin{figure}
\center{\includegraphics[width=0.99\columnwidth]{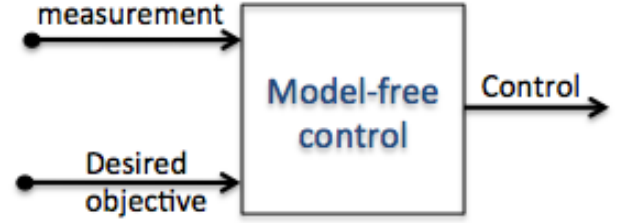}}
 \caption{General scheme of the input/output control \label{model1}}
\end{figure}
Implementing our model-free control and the related iPI controller
is straightforward:\footnote{See Fig. \ref{model1} for a
corresponding block diagram scheme representation.}
\begin{itemize}
\item The gains $K_P$ and $K_I$ are easily tuned thanks to the first
order ultra-local model \eqref{densite}.
\item The remarkable robustness properties follow from the excellent estimation $[F]_e$
of $F$.\footnote{See Section \ref{estF}.}
\item the generation of the desired trajectory (density) $\rho^\star$ is achieved
thanks to the following algorithm:
\begin{itemize}
\item Let $V_{\text{filtered}}$ be the filtered mean speed and
$V_{\text{threshold}}$ the speed threshold.\footnote{Concrete
studies (see, \textit{e.g.}, \cite{cete}) have demonstrated that the
level of service is highly degraded and the congestion phenomenon is
at its maximum, when the mean speed of individual vehicle is about
$30 \ km/h$. The threshold of discomfort is reached, when this speed
is equal to $85 \ km/h$ (see e.g. \cite{cete}).}
\item $\rho_{\text{d}0}$, $\rho_{\text{inc}}$, $\rho_{\text{dec}}$
denote respectively the initial
density,  the increment and decrement of the desired density.
\item If $V_{\text{filtered}}>V_{\text{threshold}}$, then
$\rho^\star=\rho_{\text{d}0}+\rho_{\text{inc}}$.
\item If $V_{\text{filtered}} < V_{\text{threshold}}$, then $\rho^\star=\rho_{\text{d}0}-
\rho_{\text{dec}}$.
\end{itemize}
\end{itemize}

\subsection{Simulation results}
\begin{figure*}
\center{\includegraphics[width=2.05\columnwidth]{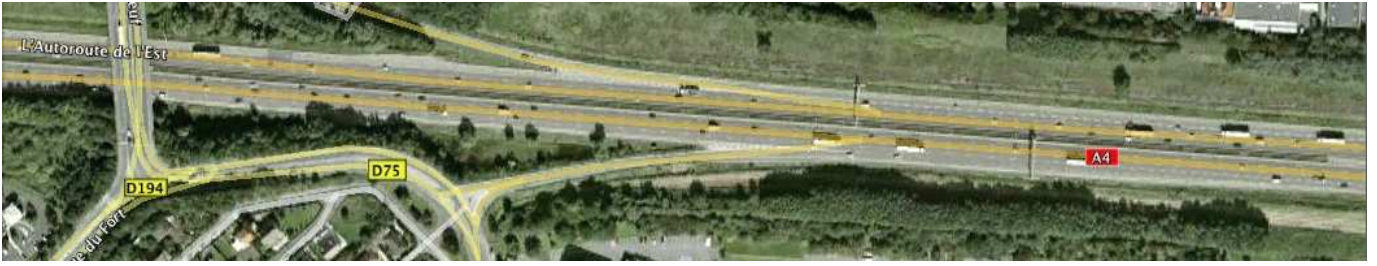}}
 \caption{Aerial picture of the studied site (Source DiRIF)\label{site}}
\end{figure*}

\begin{figure}
\center{\includegraphics[width=1.1\columnwidth]{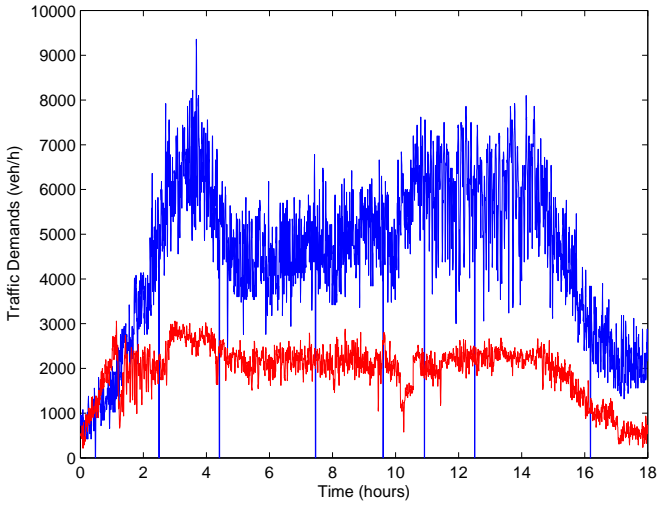}}
 \caption{Traffic demands: \textcolor{blue}{(--)} mainstream, \textcolor{red}{(--)} on-ramp\label{demands}}
\end{figure}

\begin{figure}
\center{\includegraphics[width=1.1\columnwidth]{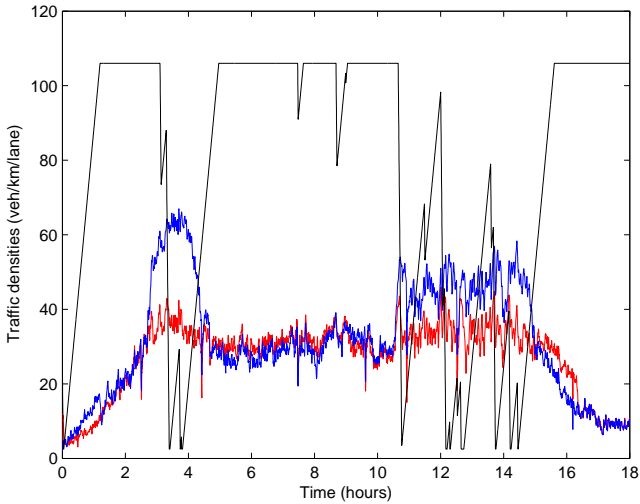}}
 \caption{Traffic densities evolutions: \textcolor{blue}{(--)} no-control case,
 \textcolor{red}{(--)} control case, {(--)} $\rho^\star$ \label{adaptive_dens}}
\end{figure}

\begin{figure}
\center{\includegraphics[width=1.1\columnwidth]{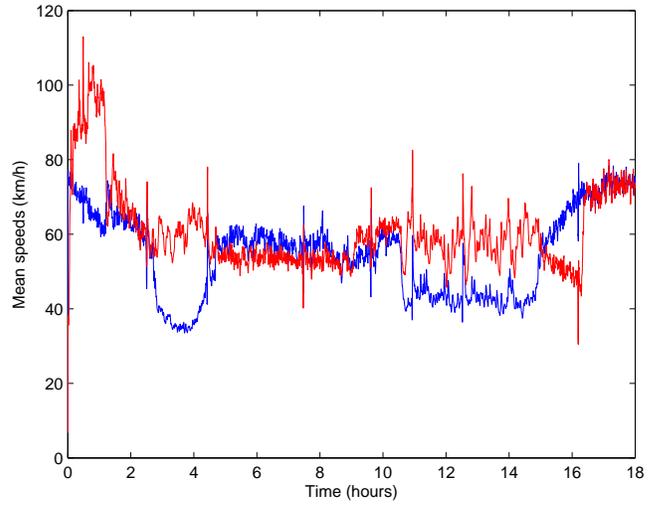}}
 \caption{Mean speeds evolutions: \textcolor{blue}{(--)} no-control case,
 \textcolor{red}{(--)} control case\label{adaptives_speeds}}
\end{figure}

\begin{figure}
\center{\includegraphics[width=1.1\columnwidth]{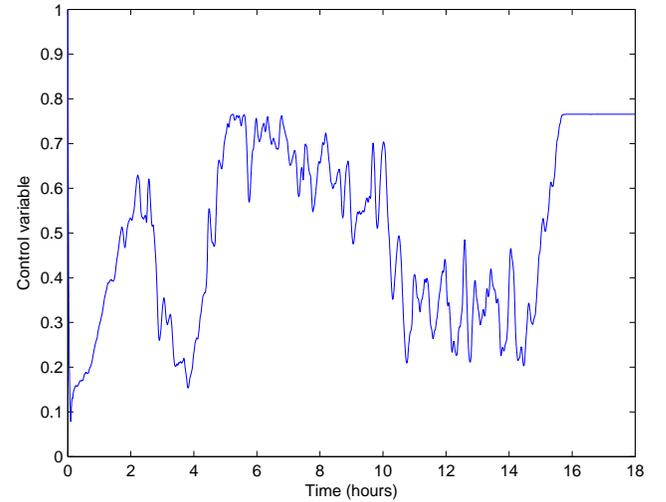}}
 \caption{Control variable\label{adap_cde}}
\end{figure}

Our computer simulations are based on numerical data which are
collected from the French freeway $A4Y$ with one on-ramp (see Fig.
\ref{site} and Fig. \ref{demands}). The software \emph{METANET}
(\cite{metanet83}) is utilized.\footnote{It is based on a second
order macroscopic model.} Although the measurements, \textit{i.e.},
the traffic volume in $veh/h$, are quite poor and noisy, the
performances of our iPI controller (Fig. \ref{adap_cde}) are good.
Congestions are alleviated as soon as they appear (Fig.
\ref{adaptive_dens} and Fig. \ref{adaptives_speeds}).

%In order to take into account the congestion periods, the numerical simulations were conducted between $5 \ am$ and $ 11 \ pm$.

% Fig.~\ref{demands} depicts the traffic demands at the two origins (mainstream section and on-ramp).
%Notice that the traffic measurements are poor and somewhat noisy. During the rush hours (7:30 am - 10:30 am) and (3:30pm - 8:30pm), we remark severe congestions, due to an important demands.
%
%Fig.~\ref{adaptive_dens} and Fig.~\ref{adaptives_speeds} show the evolution of the traffic densities and the mean speed in both cases: no-control and model-free control. We remark that in the no-control case, the traffic flow  reaches the congested mode, which is confirmed by the drop of the mean speed. In this situation, the control algorithm operates adaptively in order to maintain the traffic state at the fluid mode, which is demonstrated by the evolution of the control variable, Fig.~\ref{adap_cde}.

\section{Traffic flow parametric estimation}\label{pe}
\subsection{Generalities}
The macroscopic models which are used for simulation purposes, are
not only heuristic but also quite sensitive to parameter variations
and uncertainties. The only available accurate physical law is the
conservation equation. All other equations (speed equation and the
fundamental diagrams,  for instance), are based on empirical
observations and coarse approximations. The main parameters such as
the critical density and the free-flow speed are moreover subject to
variations.

%The success of any control strategy is closely related to the
%ability to take into account the parametric variation of the studied
%system.

%Although the "\textit{model-free control}" approach does not need
%any precise mathematical modeling, for simulation and validation
%purposes a more realistic mathematical model is needed.

%This leads us to emphasize, the dichotomy between the complexity of
%writing an exact, or at least a more realistic,  simulation (or
%mathematical) model and the simplicity of the implementation of the
%model-free control.

%The efficiency of  freeway ramp metering, for example, depends on
%the capacity to calibrate, on-line, the free-flow speed and the
%critical density and to track any changes in the main traffic
%characteristics. Although a simple estimation algorithm permits us
%to deal with such problem, for simulation purpose, the problem
%becomes more complex.

%Let us expose such identification method using the tool described in section \ref{ns}.
\subsection{Fundamental diagram}\label{may}
\begin{figure}
\center{\includegraphics[angle=-90,width=1.01\columnwidth]{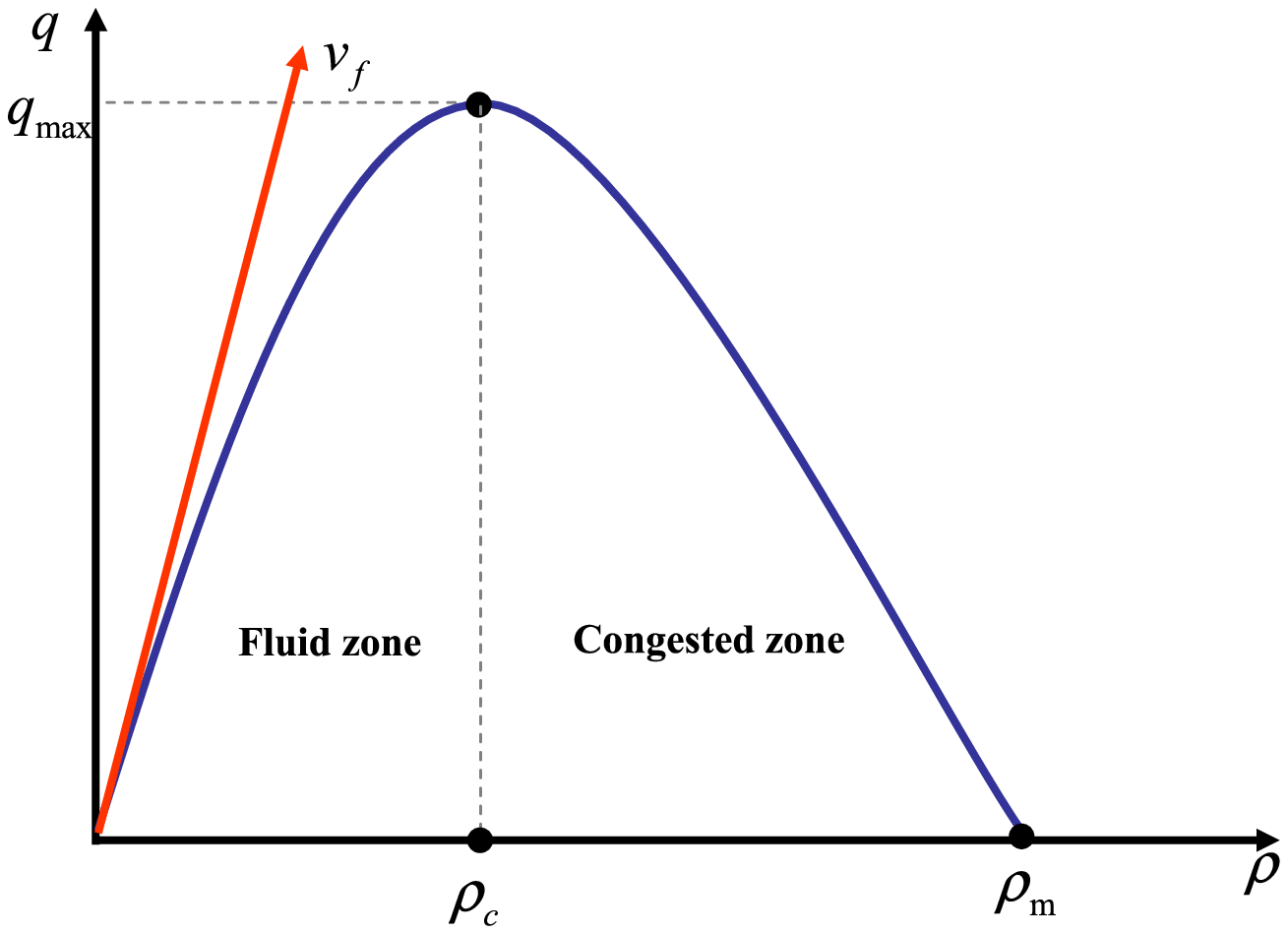}}
 \caption{Fundamental diagram example \label{diagram}}
\end{figure}

The \emph{fundamental diagram} due to \cite{May} is given by
\begin{equation}\label{df}
\boxed{V(\rho_i)=v_f
\exp\left(-\frac{1}{a}\left(\frac{\rho_i}{\rho_{c}}\right)^a
\right)}
\end{equation}
where
\begin{itemize}
\item $\rho_i$ is the density of the segment $i$,
\item $V$ is the corresponding the mean speed,
\item $v_f$ is the free-flow speed,
\item $\rho_{c}$ the critical density,
\item $a$ is a model parameter.
\end{itemize}
Although it is of heuristic nature (see Fig. \ref{diagram}), {\it
i.e.}, not derived from physical laws, it provides important
parameters for the macroscopic modeling which is used for the
numerical simulations and to identify congestion and fluid zones.

\subsection{Identification}
The existing results in signal processing (see
\cite{mexico,lobry,mboup0}) should be extended to the arbitrary
exponent $a$ in Equation \eqref{df}.
\subsubsection{New setting}
Rewrite Equation \eqref{df} in the following form
\begin{equation}\label{b}
V(\rho_i) = v_f \exp \left[ - K \rho^a \right]
\end{equation}
where $K=\frac{1}{a\rho_c^a}$. The equality
$$
\rho_{c} = \sqrt[a]{\frac{K}{a}}
$$
shows that $\rho_{c}$ may be deduced at once from $a$ and $K$.

\textbf{Notation}: If $G$ is a function of $\rho_i$, write
$G_{\rho_i}$ its derivative with respect to $\rho_i$.

Write $W$ the logarithmic derivative of $V$ with respect to
$\rho_i$:
\begin{equation}\label{K}
W = \frac{V_{\rho_i}}{V} = - K a \rho^{a - 1}
\end{equation}
Thus
\begin{equation*}\label{log}
 \frac{W_{\rho_i}}{W} =  \frac{a - 1}{\rho_i}
\end{equation*}
The identifiability of $a$ follows at once. Equation \eqref{K}
provides $K$ and Equation \eqref{b} $v_f$.
\begin{remark}
Note that the second order derivative $V_{\rho_i^2}$ of $V$  is
needed.
\end{remark}
\subsubsection{Derivation with respect to time}
Consider $\rho_i$ and, then, $V$ as functions of time $t$. The time
derivatives are obtained using the following expression:
\begin{equation}
\label{ratio}
V_{\rho_i}= \frac{\dot{V}}{\dot{\rho_i}}
\end{equation}
The numerical derivation of noisy signals, developed in
(\cite{nl,mboup}), has been already successfully implemented in many
concrete applications (see, for example in intelligent
transportation systems, \cite{wash,grenoble,vil1,vilcep}).
%\subsubsection{Derivatives of a noisy signal}\label{ns}
In order to summarize the general principles, let us start with the
first degree polynomial time function $p_1 (t) = a_0 + a_1 t$, $t
\geq 0$, $a_0, a_1 \in \mathbb{R}$. Rewrite it thanks to classic
operational calculus (see, \textit{e.g.}, \cite{yosida}) $p_1$ as
$P_1 = \frac{a_0}{s} + \frac{a_1}{s^2}$. Multiply both sides by
$s^2$:
\begin{equation}\label{1}
s^2 P_1 = a_0 s + a_1
\end{equation}
Take the derivative of both sides with respect to $s$, which
corresponds in the time domain to the multiplication by $- t$:
\begin{equation}\label{2}
s^2 \frac{d P_1}{ds} + 2s P_1 = a_0
\end{equation}
The coefficients $a_0, a_1$ are obtained via the triangular system
of equations (\ref{1})-(\ref{2}). We get rid of the time
derivatives, \textit{i.e.}, of $s P_1$, $s^2 P_1$, and $s^2 \frac{d
P_1}{ds}$, by multiplying both sides of Equations
(\ref{1})-(\ref{2}) by $s^{ - n}$, $n \geq 2$. The corresponding
iterated time integrals are low pass filters which attenuate the
corrupting noises, which are viewed as highly fluctuating phenomena
(\cite{bruit}). A quite short time window is sufficient for
obtaining accurate values of $a_0$, $a_1$.

The extension to polynomial functions of higher degree is
straightforward. For derivatives estimates up to some finite order
of a given smooth function $f: [0, + \infty) \to \mathbb{R}$, take a
suitable truncated Taylor expansion around a given time instant
$t_0$, and apply the previous computations. Resetting  and utilizing
sliding time windows permit to estimate derivatives of various
orders at any sampled time instant.

\subsection{Computer experiments}\label{sim1}
\begin{figure}
\subfigure[Density evolution $\rho$ vs time in
seconds]{\rotatebox{-0}{\includegraphics[width=1.1\columnwidth]{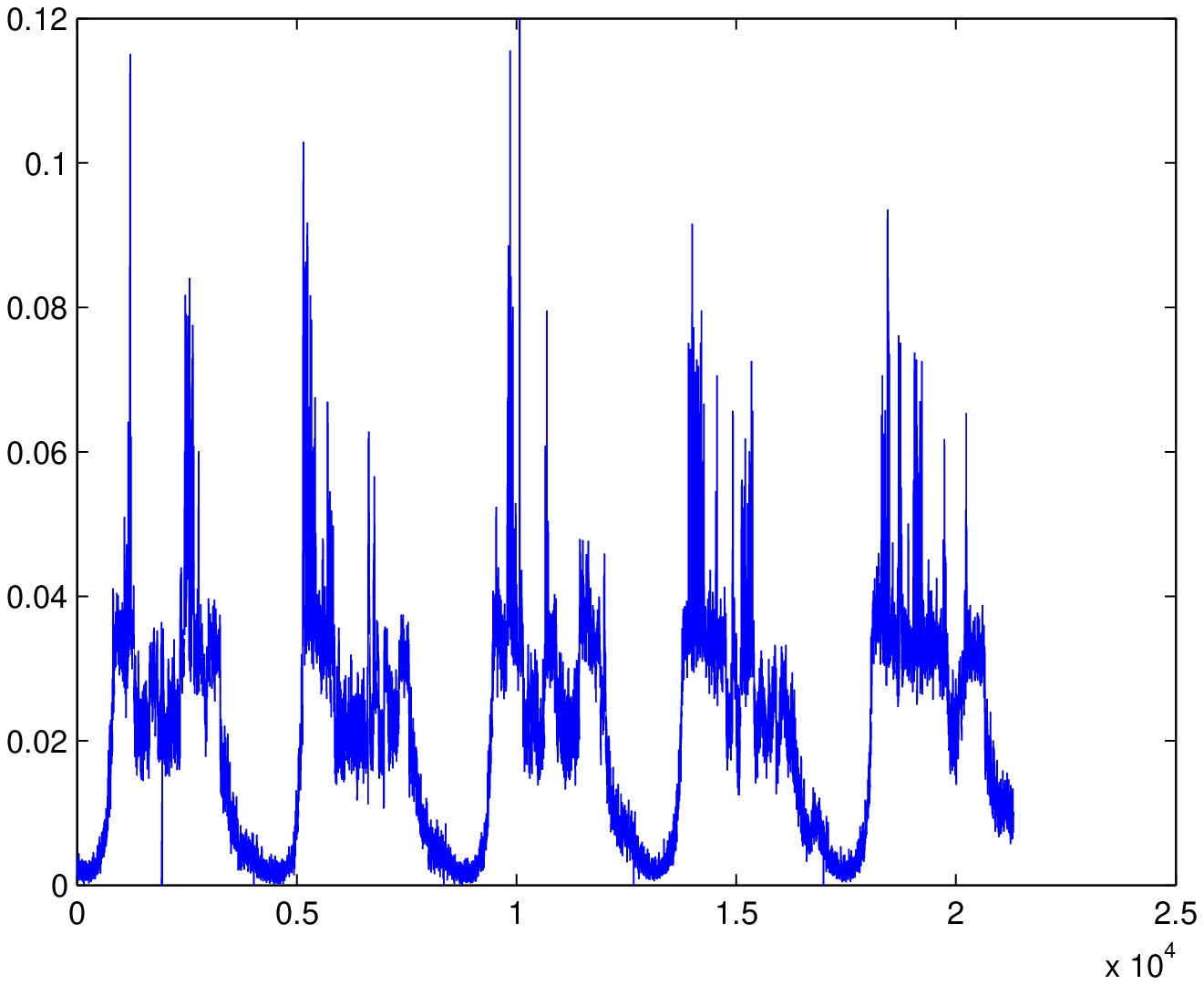}}}
\subfigure[Speed evolution $v$ vs time in
seconds]{\rotatebox{-0}{\includegraphics[width=1.1\columnwidth]{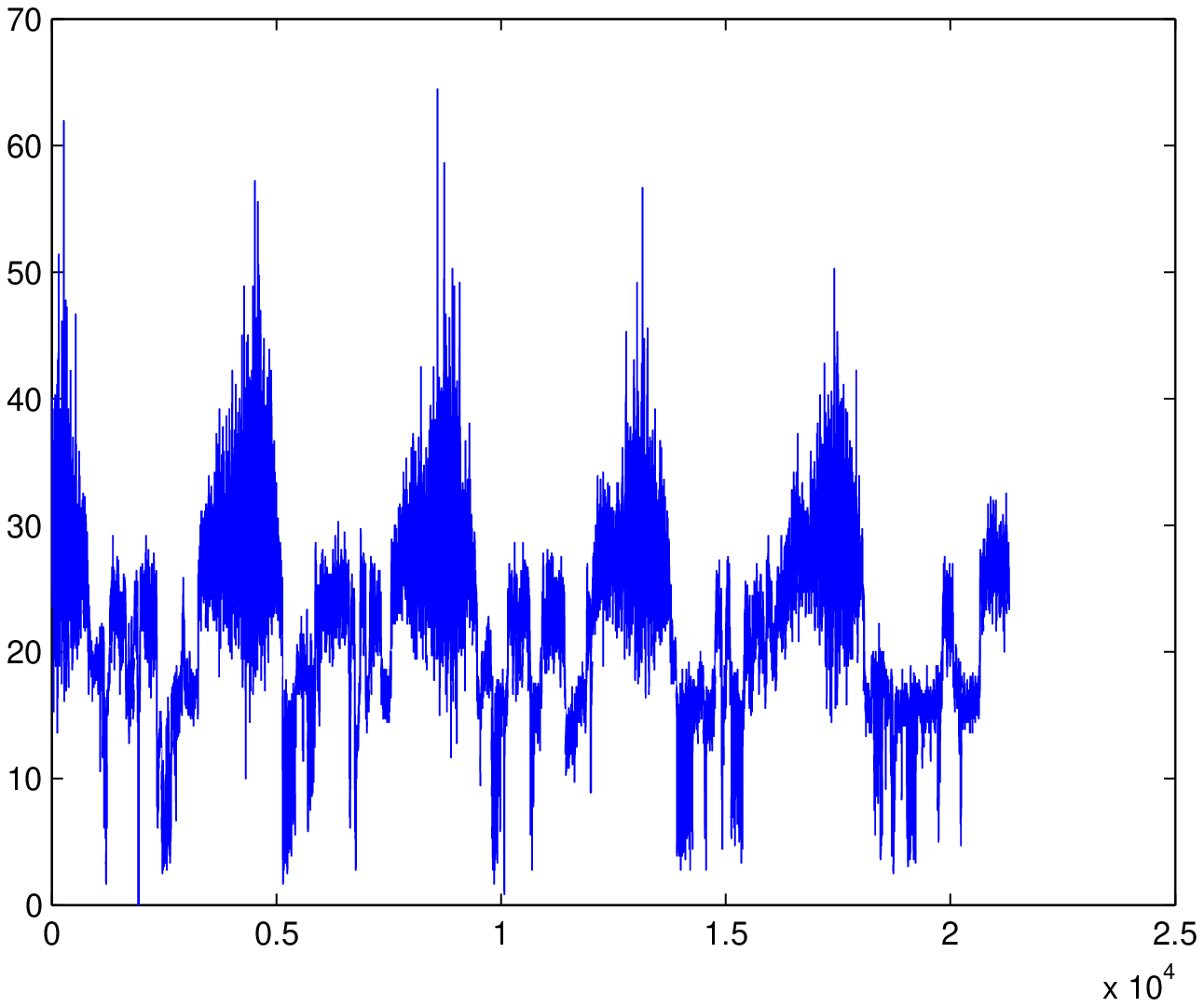}}}
 \caption{Measured variables}
 \label{mes1}
\end{figure}

\begin{figure}
\subfigure[Critical density estimation $\rho_c$ vs time in
seconds]{\rotatebox{-0}{\includegraphics[width=1.1\columnwidth]{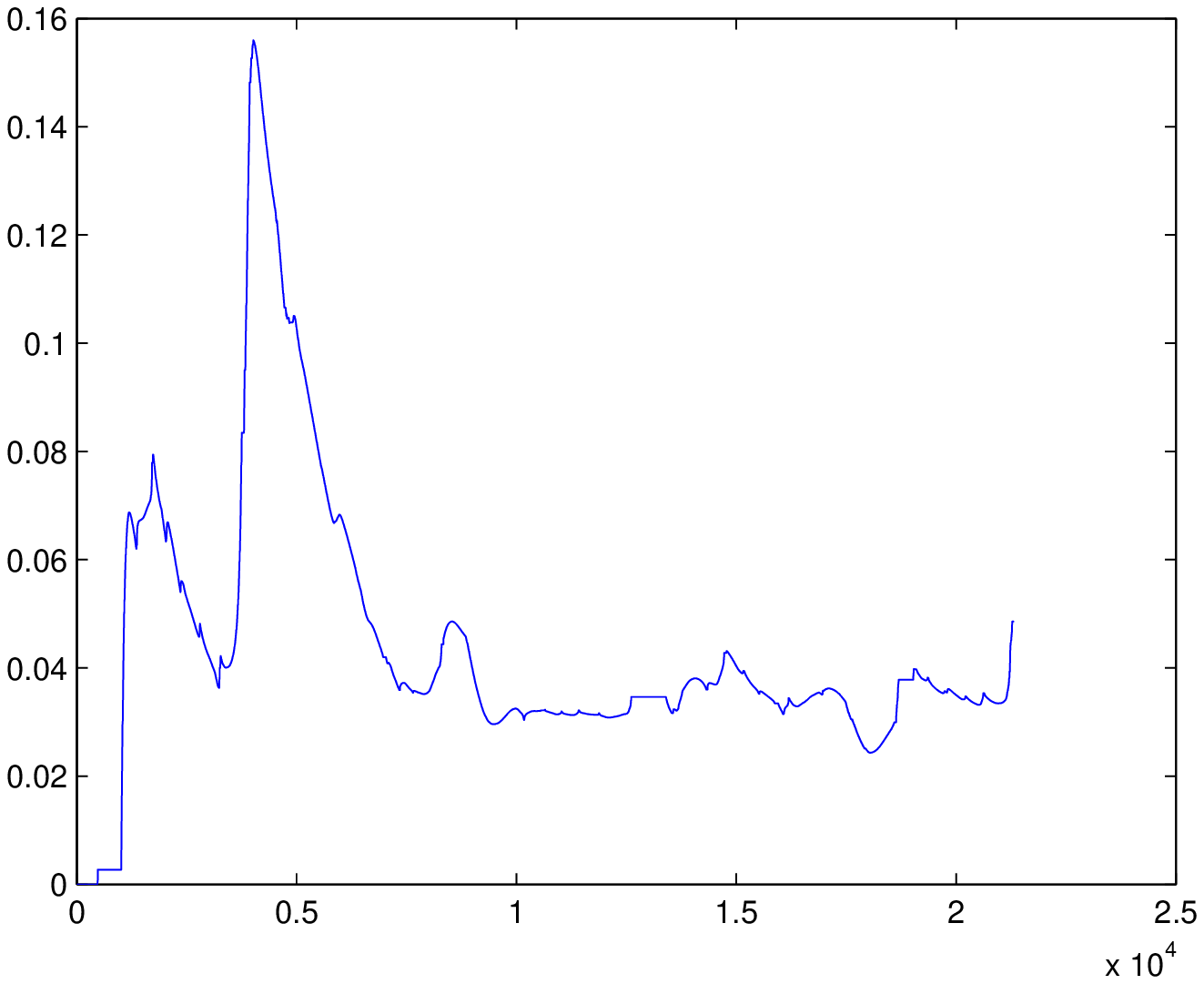}}}
\subfigure[Free-flow speed estimation $v$ vs time in
seconds]{\rotatebox{-0}{\includegraphics[width=1.1\columnwidth]{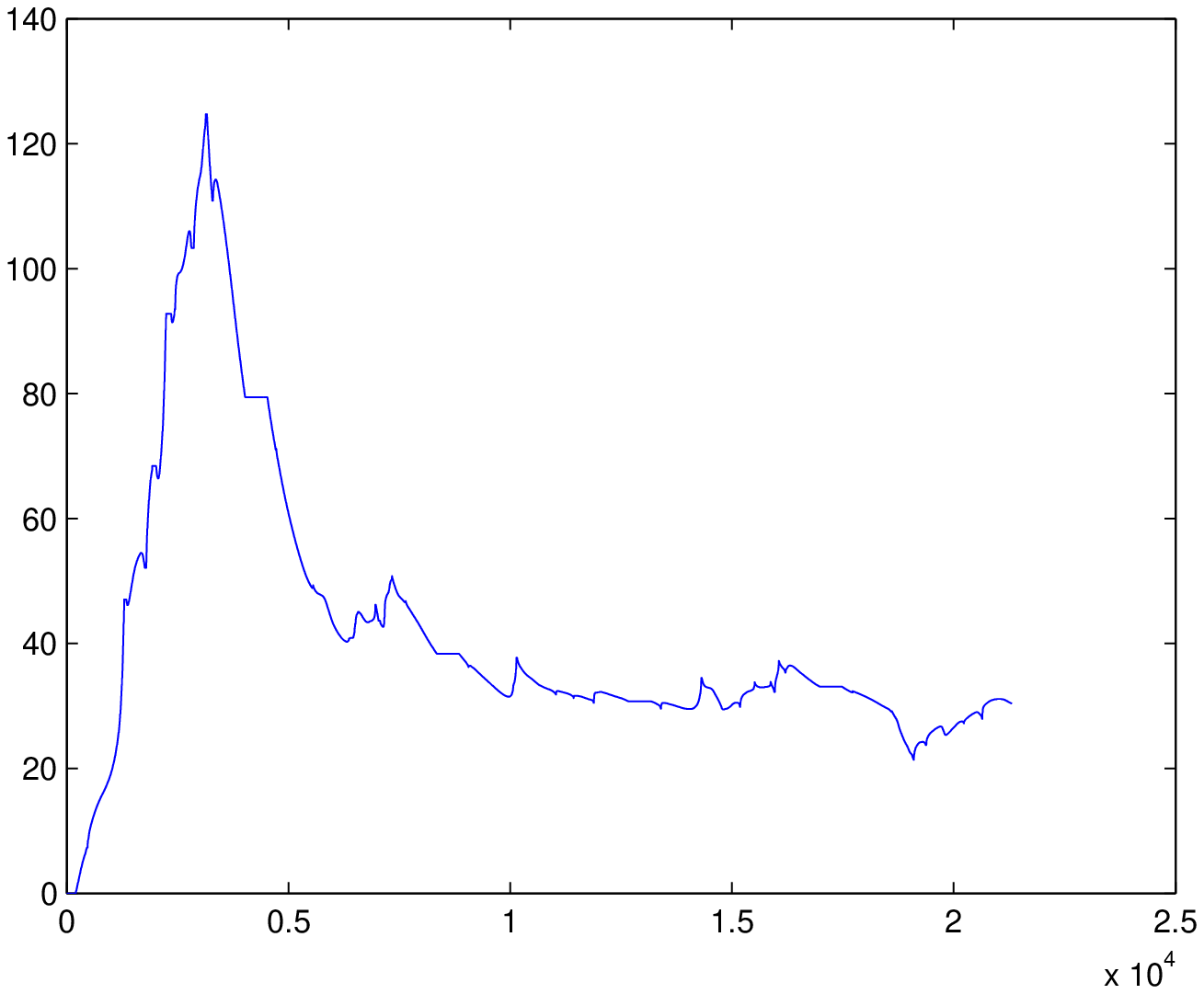}}}
\subfigure[Estimation of  $a$ vs time in
seconds]{\rotatebox{-0}{\includegraphics*[width=1.1\columnwidth]{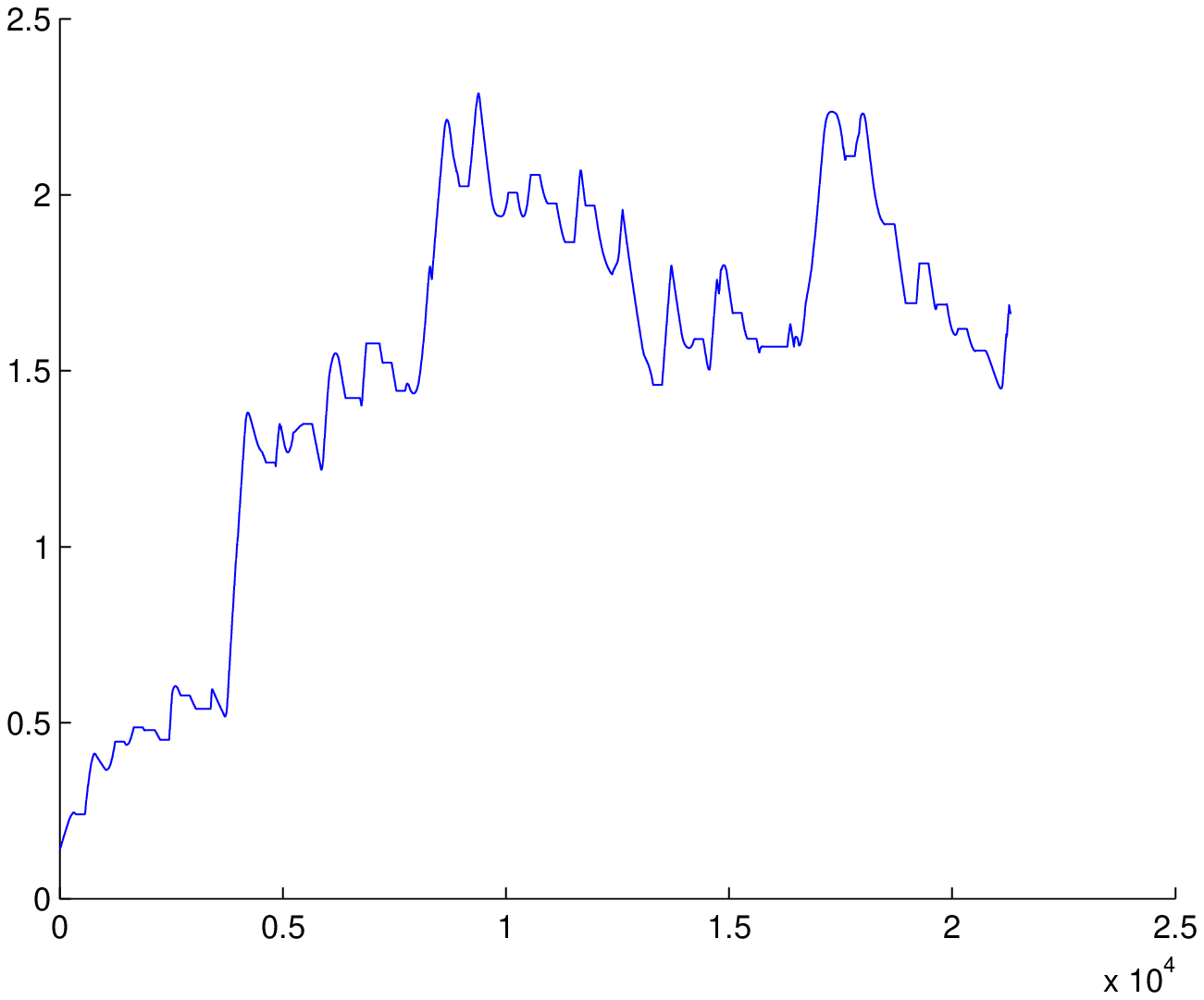}}}
 \caption{Estimated parameters}\label{estim1}
\end{figure}

The measurements concern the evolution of the mean speed and of the
traffic occupancy depicted in Fig. \ref{mes1}.\footnote{The
occupancy measurements are transformed into traffic density for
simulations purposes.}
The data are
\begin{itemize}
\item provided during five days with a sampling period of $20$
seconds,
\item quite poor and noisy.
\end{itemize}
The modeling approximation,\footnote{As already stated in Section
\ref{may} the fundamental diagram is only heuristic.} and the
numerical singularities, which are unavoidable in such a real-time
setting, explain why our estimates do fluctuate to some extent. The
results depicted in Fig. \ref{estim1} do however show a satisfactory
``practical'' convergence towards values which are suitable for our
simulation purposes.
%shows that, despite the successive losses of
%identifiability, it establishes a remarkable convergence of the
%estimations toward the approached values of the critical density,
%Fig.~\ref{estim1}-(a), the free-flow speed, $v_f$
%Fig.~\ref{estim1}-(b) and the model parameter $a$,
%Fig.~\ref{estim1}-(c).
%\section{Model-free ramp metering}\label{mfrm}
%\begin{table}[h]
%\caption{\label{tab1} Simulation parameters}
%\begin{center}
%\begin{tabular}{|c|c|}
%\hline $\rho _{c}$ & $29.5 \ \text{veh/km/lane} $ \\ \hline $\rho
%_{max}$ & $180 \ \text{veh/km/lane} $ \\ \hline $v _{f}$ & $100.1 \ \text{km/h} $ \\
%\hline $a$ & $2.997 $ \\ \hline
%\end{tabular}
%\vspace{-0.5cm}
%\end{center}
%\end{table}
%
%where $\rho_{max}$ represents the maximum density.
%

\section{Conclusion}\label{conclusion}
The control design described in this paper will soon be implemented
in practice. The pending patent\footnote{See the acknowledgement
below.} prevents us unfortunately from discussing future
developments in traffic control.

From a theoretical standpoint two major points might however be
stressed:
\begin{enumerate}
\item The control of freeway ramp metering and of hydroelectric power plants
is approached almost exclusively in the existing academic literature
via a rather complex modeling where partial differential equations
are often utilized.\footnote{See, \textit{e.g.}, some references in
\cite{agadir}.} The present work confirms what has already been
obtained for hydroelectric power plants by \cite{edf}, namely that
an elementary model-free control is enough for obtaining excellent
results.\footnote{See the conclusion in (\cite{marseille}) for a
more thorough analysis.}
\item Computer simulations do necessitate a quite realistic modeling
which implies a more subtle mathematical setting than model-free
control. This fact which is stressed here for the first time will be
further studied in the future. It might lead to a profound
epistemological revolution in applied mathematics and, more
generally, in applied sciences, the consequences of which are not
yet clear.
\end{enumerate}

{\begin{ack}\label{ack} Work partially supported by a research
contract, entitled \textit{Application de la commande sans mod\`{e}le au
contr\^{o}le d'acc\`{e}s}, between the ``Direction des Routes
d'\^{I}le-de-France'' (DiRIF) and the ``\'Ecole polytechnique''. The
``\'Ecole polytechnique'', the ``Universit\'{e} d'Artois'', the
``Universit\'{e} Henri Poincar\'{e}'' (Nancy 1) and the ``Centre National de
la Recherche Scientifique'' (CNRS) took a patent, which is pending
(n$^\circ$ FR1151604), on the regulation system described in this
paper.
\end{ack}}

%\bibliography{ifacconf}             % bib file to produce the bibliography

\begin{thebibliography}{4}
\providecommand{\natexlab}[1]{#1}
\providecommand{\url}[1]{\texttt{#1}}
\providecommand{\urlprefix}{URL }
\expandafter\ifx\csname urlstyle\endcsname\relax
  \providecommand{\doi}[1]{doi:\discretionary{}{}{}#1}\else
  \providecommand{\doi}{doi:\discretionary{}{}{}\begingroup
  \urlstyle{rm}\Url}\fi

\bibitem[{Aboua\"{\i}ssa, Fliess, Iordanova \& Join(2011a)}]{douai}
Aboua\"{\i}ssa, H., Fliess, M., Iordanova, V., Join, C. (2011a).
\newblock Vers une caract\'{e}risation non lin\'{e}aire d'un r\'{e}seau
autoroutier. \newblock \emph{3$^{es}$ J. Identif. Mod\'elisation
Exp\'erimentale}, Douai.
\newline
({\small \tt http://hal.archives-ouvertes.fr/hal-00572818/en/})

\bibitem[{Aboua\"{\i}ssa, Fliess, Iordanova \& Join(2011b)}]{agadir}
Aboua\"{\i}ssa, H., Fliess, M., Iordanova, V., Join, C. (2011b).
\newblock Prol\'{e}gom\`{e}nes \`{a} une r\'{e}gulation sans mod\`{e}le du trafic
autoroutier.
\newblock \emph{ Conf. m\'{e}dit. ing\'{e}nierie s\^{u}re syst\`{e}mes complexes},
Agadir.
\newline ({\small \tt http://hal.archives-ouvertes.fr/hal-00585442/en/})

\bibitem[{Alvarez, Horowitz \& Li(1999)}]{alvarez}
Alvarez, L., Horowitz, R., P. Li, P. (1999). Traffic flow control in
automated highway systems. \emph{Control Engin. Practice}, 7,
1071--1078.

\bibitem[{Andary, Chemori \& Benoit(2012)}]{acc}
Andary, S., Chemori, A., Benoit, M. (2012).
\newblock A dual model-free control of underactuated mechanical
systems -- Application to the inertia wheel inverted pendulum with
real-time experiments.
\newblock \emph{Amer. Control Conf.}, Montr\'{e}al.

\bibitem[{d'Andr\'{e}a-Novel, Boussard, Fliess, el Hamzaoui, Mounier \& Steux(2010)}]{cifa10}
d'Andr\'{e}a-Novel, B.,  Boussard, C.,  Fliess, M., el Hamzaoui, O.,
Mounier, H., Steux, B. (2010). \newblock Commande sans mod\`{e}le de
vitesse longitudinale d'un v\'{e}hicule \'{e}lectrique. \newblock {\em
6$^{e}$ Conf. Internat. Francoph. Automatique}, Nancy.
\newline ({\small \tt
http://hal.archives-ouvertes.fr/inria-00463865/en/})

\bibitem[{d'Andr\'{e}a-Novel, Fliess, Join, Mounier \& Steux(2010)}]{pid}
d'Andr\'{e}a-Novel, B., Fliess, M., Join, C., Mounier, H., Steux, B.
(2010). \newblock A mathematical explanation via ``intelligent'' PID
controllers of the strange ubiquity of PIDs. \newblock
\emph{18$^{th}$ Medit. Conf. Control Automat.}, Marrakech.
\newline ({\small \tt
http://hal.archives-ouvertes.fr/inria-00480293/en/}).


\bibitem[{Bellemans, De Schutter \& De Moor(2006)}]{bellemans}
Bellemans, T., De Schutter, B., De Moor, B. (2006). \newblock Model
predictive control for ramp metering of motorway traffic: A case
study. \newblock \textit{Control Engin. Practice}, 14, 757-767.

\bibitem[{Blandin, Work, Goatin, Piccoli \& Bayen(2011)}]{bayen}
Blandin, S., Work, D., Goatin, P., Piccoli, B., Bayen, A. (2011)
\newblock A general phase transition model for vehicular traffic.
\newblock \emph{SIAM J. Appl. Math.}, 71. 107--127.

\bibitem[{Boukhnifer \& Haj-Salem(2010)}]{habib_cifa}
Boukhnifer, M., H. Haj-Salem, H. (2010). \newblock Evaluation
op\'erationnelle de la r\'egulation d'acc\`es sur les autoroutes de
l'IDF bas\'ee sur la strat\'egie ALINEA. \newblock \emph{6$^e$ Conf.
Internat. Francoph. Automat.}, Nancy, 2010.

\bibitem[{Chiang \& Juang(2008)}]{chiang}Chiang Y.-H., Juang, J.-C. (2008).
\newblock Control of
freeway traffic flow in unstable phase by $H_\infty$ theory.
\emph{IEEE Trans. Intel. Transport. Systems}, 9, 193--208.

\bibitem[{Cete M\'editerran\'ee - Les \'etudes (2006)}]{cete}
CETE M\'edirann\'ee (2006).
\newblock Niveaux de service des r\'eseaux routiers en PACA : \'Etat en 2004,
\'Evolution estim\'ee \`a l'horizon 2020.
\newblock \emph{Rapport d'\'Etudes}.

\bibitem[{Choi, d'Andr\'{e}a-Novel, Fliess, Mounier \& Villagra(2009)}]{buda}
Choi, S., d'Andr\'{e}a-Novel, B., Fliess, M., Mounier, H., Villagra, J.
(2009).
\newblock Model-free control of automotive engine and brake
for Stop-and-Go scenarios.
\newblock {\em 10$^{th}$ Europ. Control
Conf.}, Budapest. \newline ({\small \tt
http://hal.archives-ouvertes.fr/inria-00395393/en/})

\bibitem[{De Miras, Riachy, Fliess, Join \& Bonnet(2012)}]{cifa12}
De Miras, J., Riachy, S., Fliess, M., Join, C., Bonnet, S. (2012).
\newblock Vers une commande sans mod\`ele d'un palier magn\'etique.
\newblock \emph{7$^e$ Conf. Internat. Francoph. Automatique},
Grenoble. \newline ({\small \tt
http://hal.archives-ouvertes.fr/hal-00682762/en/})


\bibitem[{Fliess(2006)}]{bruit}Fliess, M. (2006).
\newblock Analyse non standard du bruit.
\newblock \emph{C.R. Acad. Sci. Paris Ser. I}, 342, 797--802.

\bibitem[{Fliess(2008)}]{lobry}Fliess, M. (2008). \newblock Critique du rapport signal \`{a} bruit
en communications num\'{e}riques. \newblock \emph{ARIMA}, 9, 419-429.
\newline ({\small \tt
http://hal.archives-ouvertes.fr/inria-00311719/en/})


\bibitem[{Fliess \& Join(2008)}]{esta}
Fliess, M., Join, C. (2008).
\newblock Commande sans mod\`{e}le et
commande \`{a} mod\`{e}le restreint.
\newblock \emph{e-STA}, 5
(n$^\circ$ 4), 1--23. \newline ({\small \tt
http://hal.archives-ouvertes.fr/inria-00288107/en/})

\bibitem[{Fliess \& Join(2009)}]{malo}Fliess, M., Join, C. (2009).
\newblock  Model-free control
and intelligent PID controllers: towards a possible trivialization
of nonlinear control? \newblock\emph{15$^{th}$ IFAC Symp. System
Identif.}, Saint-Malo. \newline ({\small \tt
http://hal.archives-ouvertes.fr/inria-00372325/en/})

\bibitem[{Fliess, Join \& Riachy(2011)}]{marseille}Fliess, M., Join, C., Riachy S. (2011).
\newblock Rien de plus utile qu'une bonne th\'eorie: la commande sans
mod\`ele.
\newblock {\em JD-JN MACS.}, Marseille. \newline
({\small \tt http://hal.archives-ouvertes.fr/hal-00581109/en/})


\bibitem[{Fliess, Join \& Sira-Ram\'{\i}rez(2008)}]{nl}
Fliess, M., Join, C., Sira-Ram\'{\i}rez, H. (2008).
\newblock Non-linear estimation is easy.
\newblock \emph{Int. J. Model. Identif. Control}, 4, 12--27.
\newline ({\small \tt
http://hal.archives-ouvertes.fr/inria-00158855/fr/})


\bibitem[{Fliess, Mboup, Mounier \& Sira-Ram\'{\i}rez(2003)}]{mexico}
Fliess, M., Mboup, M., Mounier, H., Sira-Ram\'{\i}rez, H. (2003).
\newblock Questioning some paradigms of signal processing via concrete
examples. \newblock In H. Sira-Ram\'{\i}rez, G. Silva-Navarro (Eds):
\emph{Algebraic Methods in Flatness, Signal Processing and State
Estimation}, Editorial Lagares, pp. 1-21. \newline ({\small \tt
http://hal.archives-ouvertes.fr/inria-00001059/en/})

\bibitem[{Formentin, de Filippi,
Tanelli \& Savaresi(2010)}]{nolcos}Formentin, S., de Filippi, P.,
Tanelli, M., Savaresi, S. (2010). \newblock Model-free control for
active braking systems in sport motorcycles. {\em 8$^{th}$ IFAC
Symp. Nonlinear Control Systems}, Bologne, 2010.

\bibitem[{G\'{e}douin, Delaleau, Bourgeot,
Join, Arab-Chirani \& Calloch(2011)}]{brest}G\'{e}douin, P.-A.,
Delaleau, E., Bourgeot, J.-M., Join, C., Arab-Chirani, S., Calloch
S. (2011).
\newblock Experimental comparison of classical pid and model-free
control: position control of a shape memory alloy active spring.
\newblock {\em Control Engin. Practice}, 19, 433--441.

\bibitem[{Ghosh \& Li(2010)}]{gl}Ghosh, S., Li T.S.
(2010).
\newblock \emph{Intelligent Transportation Systems: Smart and Green
Infrastructure Design}.
\newblock CRC Press.

\bibitem[{Hadj-Salem, Blosseville, Dav\'{e}e \& Papageorgiou(1988)}]{hist1}
Hadj-Salem, H., Blosseville, J.-M., Dav\'{e}e, M.M., Papageorgiou, M.
(1988). \newblock ALINEA: Un outil de r\'{e}gulation d'acc\`{e}s isol\'{e} sur
autoroute -- \'Etude comparative sur site r\'{e}el. \newblock
\emph{Rapport INRETS n$^\circ$ 80}, Arcueil.


\bibitem[{Hadj-Salem, Blosseville \& Papageorgiou(1990)}]{hist2}
Hadj-Salem, H., Blosseville, J.-M., Papageorgiou, M. (1990).
\newblock ALINEA - a local feedback control law for on-ramp metering: a
real life study. \newblock \emph{3$^{rd}$ IEE Intern. Conf. Road
Traffic Control}, London, pp. 194--198.

\bibitem[{Hegyi, De Schutter \& Hellendoor}(2005)]{Hegyi_r}
Hegyi, A., De Schutter,  B.,  Hellendoor, H. (2005).
\newblock Model predictive control for optimal coordination of
ramp metering and variable speed limits.
\newblock \emph{Transportation Research C}, 13, 185--209.


\bibitem[{Join, Masse \& Fliess(2007)}]{psa}Join, C., Masse, J., Fliess, M. (2008).
\newblock \'Etude pr\'{e}liminaire d'une commande sans
mod\`{e}le pour papillon de moteur. \newblock \emph{J. europ. syst.
automat.}, 42, 337--354. \newline ({\small \tt
http://hal.archives-ouvertes.fr/inria-00187327/en/})

\bibitem[{Join, Robert \& Fliess(2010)}]{edf}
Join, C., Robert, G., Fliess, M. (2010). \newblock Vers une commande
sans mod\`{e}le pour am\'{e}nagements hydro\'{e}lectriques en cascade. \newblock
\emph{6$^e$ Conf. Internat. Francoph. Automat.}, Nancy.
\newline ({\small \tt
http://hal.archives-ouvertes.fr/inria-00460912/en/})

\bibitem[{Kachroo \& Osbay(2003)}]{ko}Kachroo, P.,
Ozbay, K. (2003).
\newblock \emph{Feedback Ramp Metering in Intelligent Transportation
Systems}. \newblock Springer.

\bibitem[{Ozbay, Yasar \& Kachroo (2006)}]{Kachroo_estim}
Ozbay, K. Yasar, I., Kachroo, P. (2006).
\newblock Improved Online Estimation Methods for a Feedback-Based Freeway Ramp Metering Strategy.
\newblock \emph{Proc. IEEE Intelligent Transportation Systems Conf.}, September, 17-20, Toronto, Canada.





\bibitem[{Kostialos, Papageorgiou \& Middelham (2001)}]{kostialos1}
Kostialos, A.,  Papageorgiou, M., Middelham, F. (2001).
\newblock Optimal coordinated ramp metering with advanced motorway optimal control.
\newblock \emph{Proc. 80$^{th}$ Annual Meeting of the Transportation Research Board}, N$^\circ$01-3125, Washington, D.C.

\bibitem[{Mammar(2007)}]{mammar}Mammar, S. (Ed.) (2007).
\newblock \emph{Syst\`{e}mes de transport intelligents: Mod\'{e}lisation,
information et contr\^{o}le}.
\newblock Herm\`{e}s-Lavoisier.

\bibitem[{May(1990)}]{May}
May, A.D. (1990). \newblock \emph{Traffic Flow Fundamentals}.
\newblock Prentice-Hall.

\bibitem[{Mboup(2009)}]{mboup0}Mboup, M. (2009)
\newblock Parameter estimation for signals described by differential
equations. \newblock \emph{Applicable Anal.}, 88, 29--52.

\bibitem[{Mboup, Join \& Fliess(2009)}]{mboup}
Mboup, M., Join, C., Fliess, M. (2009).
\newblock Numerical differentiation with annihilators in noisy environment.
\newblock \emph{Numer. Algor.}, 50, 439--467.

\bibitem[{Menhour, d'Andr\'{e}a-Novel, Boussard, Fliess \& Mounier(2011)}]{wash}
Menhour, L., d'Andr\'{e}a-Novel, B., Boussard, C., Fliess, M., Mounier,
H. (2011). \newblock Algebraic nonlinear estimation and
flatness-based lateral/longitudinal control for automotive vehicles.
\newblock \emph{14$^{th}$ Int. IEEE Conf. Intelligent Transport. Systems},
Washington. \newline ({\small \tt
http://hal.archives-ouvertes.fr/hal-00611950/en/})

\bibitem[{Menhour, d'Andr\'{e}a-Novel, Fliess \& Mounier(2012)}]{grenoble}
Menhour, L., d'Andr\'{e}a-Novel, B., Fliess, M., Mounier, H. (2012).
\newblock Commande coupl\'{e}e longitudinale/lat\'{e}rale de v\'{e}hicules par platitude
et estimation alg\'{e}brique. \newblock \emph{7$^e$ Conf. Internat.
Francoph. Automatique}, Grenoble. \newline ({\small \tt
http://hal.archives-ouvertes.fr/hal-00686653/en/})


\bibitem[{Mihaylova, Boel \& Hegyi(2009)}]{ekf}Mihaylova, L., Boel, R., Hegyi, A. (2007).
\newblock Freeway traffic estimation within recursive Bayesian framework.
\newblock \emph{Automatica}, 43, 290--300.

\bibitem[{Michel, Join, Fliess, Sicard \& Ch\'{e}riti(2010)}]{michel} Michel, L., Join, C., Fliess, M.,
Sicard, P., Ch\'{e}riti, A. (2010). \newblock Model-free control of
dc/dc converters. \newblock \emph{IEEE Compel}, Boulder.
\newline ({\small \tt http://hal.archives-ouvertes.fr/inria-00495776/en/})


\bibitem[{Papageorgiou(1983)}]{metanet83}
Papageorgiou M. (1983).
\newblock \emph{Applications of automatic control
concepts to traffic flow modeling and control}.
\newblock Springer.


\bibitem[{Papageorgiou, Diakaki, Dinopoulou, Kostialos \& Wang(2003)}]{hist4}
Papageorgiou, M., Diakaki, C., Dinopoulou, D., Kostialos, A., Wang,
Y. (2003). \newblock Review of road traffic control strategies.
\newblock \emph{IEEE Trans. Intelligent Transport. Syst.}, 91,
2043--2067.

\bibitem[{Papageorgiou, Blosseville \& Hadj-Salem(1990)}]{papageorgiou}
Papageorgiou, M., Blosseville, J.-M.,  Hadj-Salem, H. (1990).
\newblock Modelling and real-time control of traffic flow on the southern part
of boulevard peripherique in Paris -- Part I: Modelling.
\newblock \emph{Transport. Research A}, 24, 345--359.

\bibitem[{Papageorgiou, Hadj-Salem \& Blosseville(1991)}]{hist3}
Papageorgiou, M., Hadj-Salem, H., Blosseville, J-M. (1991).
\newblock ALINEA: A local feedback control law for on-ramp metering.
\newblock
\emph{Transp. Res. Record}, n$^\circ$ 1320, 58--64.

%Papageorgiou, M., Haj-Salem, H., Middelham, F.: ALINEA local ramp
%metering: Summary of field results. Transportation Research Record
%No. 1603 (1998), pp. 90-98.
\bibitem[{Papageorgiou (1998)}]{papageorgiou98}
Papageorgiou, M. (1998).
\newblock Some remarks on macroscopic traffic flow modelling.
\newblock \emph{Transport. Research}, 32, 323--329.


\bibitem[{Payne(1971)}]{Payne}
Payne,  H. J. (1971) \newblock  Models of traffic and control.
\emph{Simulation Council Proc. Math. Models Public Syst.}, vol. 1,
chap. 6, pp. 51--61.


\bibitem[{Smaragdis \& Papageorgiou(2004)}]{hist5}Smaragdis, E.,
Papageorgiou, M. (2004). \newblock A series of new local ramp
metering strategies. \newblock \emph{Transport. Research B},  38,
251--270.


\bibitem[{Villagra, d'Andr\'{e}a-Novel, Fliess \& Mounier(2009)}]{vil1}Villagra, J.,
d'Andr\'{e}a-Novel, B., Fliess, M., Mounier, H. (2009).
\newblock Robust stop-and-go control strategy: an algebraic approach
for nonlinear estimation and control. \newblock {\em Int. J. Vehicle
Autonomous Systems}, 7, 270--291. \newline ({\small \tt
http$:$//hal.archives-ouvertes.fr/inria-00419445/en/})

\bibitem[{Villagra, d'Andr\'{e}a-Novel, Fliess \& Mounier(2011)}]{vilcep}Villagra, J.,
d'Andr\'{e}a-Novel, B., Fliess, M., Mounier, H. (2011). A
diagnosis-based approach for tire-road forces and maximum friction
estimation. \newblock \emph{Control Engin. Practice}, 19, 174--184.

\bibitem[{Villagra \& Balaguer(2011)}]{vil2}Villagra, J.,
Balaguer, C. (2011). \newblock A model-free approach for accurate
joint motion control in humanoid locomotion. \newblock \emph{Int. J.
Humanoid Robotics}, 8, 27--46.

\bibitem[{Wang, Mounier, Cela \& Niculescu(2011)}]{mounier}
Wang, J., Mounier, H., Cela, A., Niculescu, S.-I. (2011).
\newblock Event driven intelligent PID controllers with applications to motion
control.
\newblock \emph{18$^{th}$ IFAC World Congress}, Milan.

\bibitem[{Wang \& Papageorgiou(2005)}]{Wang05}
Wang, Y., Papageorgiou, M. (2005).
\newblock Real-time freeway
traffic state estimation based on extended Kalman filters: a general
approach.
\newblock \emph{Transport. Research B}, 39,
141--167.

\bibitem[{Wang, Papageorgiou \& Messmer(2008)}]{Wang08}
Wang, Y., Papageorgiou, M., Messmer, A. (2008). \newblock Real-time
freeway traffic state estimation based on extended Kalman filter:
Adaptive capabilities and real data testing.
\newblock \emph{Transport. Research A}, 42, 1340--1358.

\bibitem[{Yosida(1984)}]{yosida} Yosida, K. (1984).
\newblock {\em Operational Calculus: A Theory of Hyperfunctions}
(translated from the Japanese). \newblock Springer.

\end{thebibliography}
                                                     % with bibtex (preferred)

\end{document}